\documentclass[12pt]{amsart}
\usepackage{amscd,amsmath,amssymb,amsthm,enumerate,times,ifthen}
\usepackage{epsfig,epic,eepic,mathrsfs}
\usepackage{graphicx}
\usepackage[utf8]{inputenc}
\newtheorem{theorem}{Theorem}
\newtheorem{Theorem}{Theorem}
\usepackage[T1]{fontenc}
\usepackage{tikz}
\def\Thm#1#2{\begin{theorem}\label{T#1}#2\end{theorem}}

\def\thm#1{Theorem~\ref{T#1}}

\newtheorem*{CP}{Contraction Principle}
\def\Cp#1{\begin{CP}#1\end{CP}}

\newtheorem{Lemma}{Lemma}

\newtheorem{Prposition}{Proposition}
\def\Prp#1#2{\begin{Prposition}\label{P#1}#2\end{Prposition}}
\def\prp#1{Proposition~\ref{P#1}}

\newtheorem{Corollary}{Corollary}
\def\Cor#1#2{\begin{Corollary}\label{C#1}#2\end{Corollary}}
\def\cor#1{Corollary~\ref{C#1}}

\newtheorem*{theoremn}{Theorem}

\newtheorem*{lemman}{Lemma}

\newtheorem*{corn}{Corollary}

\theoremstyle{definition}

\newcommand{\N}{\mathbb N}

\newcommand{\R}{\mathbb R}

\newcommand{\B}{\mathscr B}
\newcommand{\CC}{\mathscr C}
\newcommand{\F}{\mathscr F}
\newcommand{\HH}{\mathscr H}

\newcommand{\hatv}{\mathscr P}
\newcommand{\W}{\mathscr W}

\newcommand{\de}{\delta}
\newcommand{\De}{\Delta}

\newcommand{\sig}{\sigma}
\newcommand{\eps}{{\varepsilon}}

\newcommand{\rank}{\mathop{\hbox{\rm \tiny rank}}}
\newcommand{\rest}{\!\!\restriction}
\renewcommand{\phi}{\varphi}

\newcommand{\Eq}[2]{\ifthenelse{\equal{#1}{*}}
{\begin{equation*}\begin{aligned}[]#2\end{aligned}\end{equation*}}
{\begin{equation}\label{#1}\begin{aligned}[]#2\end{aligned}\end{equation}}}

\begin{document}

\date{\today}

%\begin{flushright}
%\emph{}
%\end{flushright}

\title{Applications of the Bielecki renorming technique}

\author[M. Bessenyei]{Mih\'aly Bessenyei}
\author[Zs. P\'ales]{Zsolt P\'ales}

\address{Institute of Mathematics,
University of Debrecen, H-4002 Debrecen, Pf.\ 400, Hungary}

\email{besse@science.unideb.hu}
\email{pales@science.unideb.hu}



\thanks{The research was supported by the J\'anos Bolyai Research Scholarship of the Hungarian Academy of Sciences, by the \'UNKP-19-4 New National Excellence Program of the Ministry for Innovation and Technology and by the 2019-2.1.11-T\'ET-2019-00049, EFOP-3.6.1-16-2016-00022 and EFOP-3.6.2-16-2017-00015 projects. The last two projects are co-financed by the European Union and the European Social Fund.}

\begin{abstract}
The renorming technique allows one to apply the Banach Contraction Principle for maps which are not contractions with
respect to the original metric. This method was invented by Bielecki and manifested in an extremely elegant proof of
the Global Existence and Uniqueness Theorem for ODEs. The present paper provides further extensions and applications of
Bielecki's method to problems stemming from the theory of functional analysis and functional equations. 
\end{abstract}

\maketitle

\section{Introduction}

The Banach Contraction Principle \cite{Ban22} provides a sufficient condition for the fixed point property of a self map of
complete metric space in terms of contractivity. However, important situations occur when contractivity cannot be guaranteed
whereas fixed point property is still expected. In such situations, the following idea may help: Find a metric in which the
original space remains complete and in which the original map becomes a contraction. Then, the Contraction Principle applies.
On the other hand, the fixed point property is a metric-independent, algebraic property. Thus our map must have a unique fixed
point.

Another standard trick is to verify that some iterate of the given map is a contraction and then the unique fixed point property
again follows from the Contraction Principle. In fact, this approach is less general than the remetrization technique: If, for
$T\colon X\to X$ there exists $k\in\N$ such that $T^k$ is a $q$-contraction of the metric space $(X,d)$, then $T$ is a
$\sqrt[k]{q}$-contraction of the metric space $(X,d_k)$, where
\Eq{*}{
  d_k(x,y):=d(x,y)+q^{-\frac1k}d(Tx,Ty)+\cdots+q^{-\frac{k-1}k}d(T^{k-1}x,T^{k-1}y) \qquad(x,y\in X).
}
Indeed, by the $q$-contractivity of $T^k$, we have 
\Eq{*}{
   d_k(Tx,Ty)
   &=d(Tx,Ty)+q^{-\frac1k}d(T^2x,T^2y)+\cdots+q^{-\frac{k-1}k}d(T^{k}x,T^{k}y)\\ 
   &\leq d(Tx,Ty)+q^{-\frac1k}d(T^2x,T^2y)+\cdots+q^{1-\frac{k-1}k}d(x,y)
   =q^{\frac1k}d_k(x,y).
}

The remetrization idea appears in the paper of Bielecki \cite{Bie56}, and manifests in an extremely elegant proof of the
Global Existence and Uniqueness Theorem for ODEs. In fact, this proof shows the unique solvability of an integral equation
which is equivalent to the original Cauchy problem. Comparing this integral equation to that of Volterra, one can immediately
discover their relationship. However, Volterra equations are handled quite differently: The standard approach is to show that
some iterates of the map determined by the Volterra equation is a contraction in the original norm.

Therefore the question arises: \emph{Can we prove these results in the same way?} In this paper, we give a positive answer to
this question by Bielecki's method. We are going to investigate the nonlinear integral equation
\Eq{*}{
 x(t)=f(t)+\int_{H(t)}K\bigl(t,s,x(s)\bigr)d\mu(s).
}
The unknown function $x$ belongs to the space of continuous functions $\CC(X,B)$, were $B$ is a Banach space and $X$ is a locally
compact topological space with Radon measure $\mu$. The integral is meant in the sense of Bochner. The domain of integration is
given by a relation $H\subseteq X^2$ whose properties will be clarified later.

Our main results provide existence and uniqueness theorems for the solvability of the equation above. The sufficient condition
that we need, the most important feature of the theorem, is the solvability of a \emph{homogeneous linear integral inequality}
which is connected to the Lipschitz property of the kernel function $K$. The advantage of this assumption is obvious: Finding a
solution to a homogeneous linear integral inequality is much easier then finding the (unique) solution of an inhomogeneous
nonlinear integral equation. Although this is not the aim of the paper, let us point out that our  assumption can also be checked
via standard numerical methods.

Let us point out, that several important particular cases of the above equation have been studied intensively. The monographs of
Corduneanu \cite{Cor91}, of Gripenberg, Londen, and Staffans \cite{GriLonSta90}, and of Guo, Lakshmikantham, and Liu \cite{GuoLakLiu96}
give an excellent overview of the topic. Corduneanu \cite{Cor91} presents a nice issue on the prehistory and the evolution of the
seminal works of Fredholm \cite{FreOmn} and Volterra \cite{VolOmn}.

Recent developments about integral equations basically extend the \emph{range} of the functions beyond Banach spaces to fuzzy spaces
\cite{Skr17} and so-called $L$-spaces \cite{Bab19}. However, the unknown functions are  defined only on intervals. In our setup, the
generalization concerns the \emph{domain}, as well. This has an immediate effect even to the classical cases: we can treat the Volterra-
and Fredholm-type equations with the Global Existence and Uniqueness Theorem \emph{simultaneously}.

The paper is organized as follows. As preliminaries, we collect the most important tools from set theory, measure theory and functional
analysis. The most important results in this section are an extension lemma which allows to change local fixed point properties to a
global one, and a regularity lemma, which corresponds to the continuity of the classical integral function. An alternative approach to
the Bochner integrability of continuous maps on compact domains is also presented. Finally, we introduce the spectral radius function
and enlighten its connection to an integral equation. In the next section, we present our main results with their proofs. Finally, in
the last three sections, we give several applications to Fredholm- and Volterra-type integral equations and to Presi\'c-type functional
equations. Our method allows us to present the Global Existence and Uniqueness Theorem of ODEs and of a Wawe-type Equation in a common,
unique form.

\section{Preliminaries}

Throughout in this paper, $\N$ and $\R_+$ stand for the set of positive integers and the set of positive reals,
respectively. The aim of this section is to give a brief overview of the needed theoretical background. In the
first well-known statement let us recall the basic fixed point theorem which was established by Banach \cite{Ban22}
in 1922.

\Cp{If $T$ is a self-map of a nonempty set $S$ such that $S$ can be equipped with a complete metric in which $T$ is a
contraction, then $T$ has a unique fixed point in $S$. Furthermore, for all $x_1\in S$, the sequence $(x_n)_{n\in\N}$ defined by
the Banach--Piccard iteration
\Eq{*}{
  x_{n+1}:=Tx_n
}
converges to the unique fixed point of $T$.}

In the sequel, some set-theoretical tools are presented. As usual, $B^X$ stands for all maps acting on $X$ and having
values in $B$. The restriction of $x\in B^X$ to a set $H\subseteq X$ is denoted by $x\rest_H$. Let $\F$ be a subset
of $B^X$ and let $T$ be a self-map of $\F$. For a subset $H$ of $X$, denote the set $\{x\rest_H\,:\, x\in\F\}$
by $\F_H$. We say that $T$ is \emph{restrictable} to $\F_H$ if
\Eq{*}{
 (Tx)\rest_H=(Ty)\rest_H\qquad\mbox{whenever}\qquad x\rest_H=y\rest_H.
}
In this case, $T_H(x\rest_H):=(Tx)\rest_H$ defines a function $T_H\colon\F_H\to\F_H$, which we call the \emph{natural
restriction} of $T$ to $\F_H$. Our first extension result gives a sufficient condition in order to local fixed point
properties be a global one.

\Prp{solex}{Let $\F\subseteq B^X$, let $T\colon\F\to\F$, and let $\HH\subseteq\hatv(X)$. Assume that
\begin{enumerate}[(i)]
 \item $\HH$ is a cover for $X$;
 \item $\bigcup\{H\in\HH\mid H\subset H_1\cap H_2\}=H_1\cap H_2$ for all $H_1,H_2\in\HH$;
 \item if $x\rest_H\in \F_H$ for all $H\in\HH$, then $x\in\F$;
 \item $T$ is restrictable to $\F_H$ for all $H\in\HH$;
 \item for all $H\in\HH$, the natural restriction $T_H$ has a unique fixed point $x_H\in\F_H$.
\end{enumerate}
Then $T$ has a unique fixed point in $\F$.}

\begin{proof}
Assume that $H_0,H\in\HH$ and $H\subset H_0$. If $x\in\F_{H_0}$, then there exists $u\in\F$ such that
$x=u\rest_{H_0}$. Since $T$ is restrictable both to $\F_{H_0}$ and $\F_H$,
\Eq{*}{
 (T_{H_0}x)\rest_H=((Tu)\rest_{H_0})\rest_H=(Tu)\rest_H=T_H(u\rest_H)=T_H(x\rest_H)
}
follows. In particular, if $x_{H_0}$ is the unique fixed point of $T_{H_0}$, we arrive at
\Eq{*}{
 x_{H_0}\rest_H=(T_{H_0}(x_{H_0}))\rest_H=T_H(x_{H_0}\rest_H).
}
Thus $x_{H_0}\rest_H$ is a fixed point of $T_H$ in $\F_H$. However, the fixed point of $T_H$ is unique,
yielding $x_{H_0}\rest_H=x_H$. This property enables us to define a function in the following way. If
$t\in X$, then there exists $H_0\in\HH$ such that $t\in H_0$. Then let $x(t):=x_{H_0}(t)$. The definition is
correct: If $H_1$ and $H_2$ share these properties, then there exists $H\in\HH$ such that $H\subset H_1\cap H_2$
and $t\in H$. Hence, using the previous observation,
\Eq{*}{
 x_{H_1}(t)=x_{H_1}\rest_{H}(t)=x_{H}(t)=x_{H_2}\rest_{H}(t)=x_{H_2}(t).
}

Obviously, $x\rest_H=x_H$ holds for all $H\in\HH$, and hence $x$ belongs to $\F$. Moreover, we show that $x$ is
a fixed point of $T$. Indeed, if $t\in X$ and $H\in\HH$ contains $t$, then
\Eq{*}{
 (Tx)(t)=(Tx)\rest_H(t)=T_H(x\rest_H)(t)=T_H(x_H)(t)=x_H(t)=x(t).
}
On the other hand, any fixed point $x\in\F$ of $T$ possesses $x\rest_H=x_H$. Therefore the uniqueness of $x_H$
provides the uniqueness of $x$, as well.
\end{proof}

Using relations instead of covering families can be more convenient: It turns out that some well-known properties of
relations imply the first three properties of \prp{solex}. We summarize these in the next result. As usual, any subset
$H$ of $X^2$ is termed a \emph{relation} on $X$. Recall that any relation $H$ induces a set-valued map $H(\cdot)$ via the
definition
\Eq{*}{
 H(t):=\{s\in X\mid (t,s)\in H\}.
}
A relation $H$ on a topological space $X$ is called \emph{strongly surjective}, if the induced set-valued map
generates an open cover:
\Eq{*}{
 \bigcup_{t\in X}H(t)^{\circ}=X.
}

\Prp{rel}{Let $H$ be a relation on a nonempty set $X$.
\begin{enumerate}[(i)]
 \item If $H$ is transitive, then $H(s)\subseteq H(t)$ whenever $s\in H(t)$.
 \item If $H$ is reflexive and transitive, then, for all $t_1,t_2\in X$,
  \Eq{*}{
   \bigcup\{H(t)\mid H(t)\subset H(t_1)\cap H(t_2)\}=H(t_1)\cap H(t_2).}
\end{enumerate}
Assume that $X,Y$ are topological spaces, and $H$ is a strongly surjective relation on $X$. If a map $x\colon X\to Y$
satisfies $x\rest_{H(t)}\in\CC(H(t),Y)$ for all $t\in X$, then $x\in\CC(X,Y)$.}

\begin{proof}
Assume that $H$ is transitive. Fix $s\in H(t)$ and choose $u\in H(s)$. Then, $(t,s)\in H$ and $(s,u)\in H$. By transitivity,
$(t,u)\in H$. Thus $u\in H(t)$, and hence $H(s)\subseteq H(t)$ follows.

Assume that $H$ is reflexive. Then, $t\in H(t)$ for all $t\in X$, in particular, for all $t\in H(t_1)\cap H(t_2)$. On the
other hand, the transitivity implies that $H(t)\subset H(t_1)\cap H(t_2)$ also holds for all $t\in H(t_1)\cap H(t_2)$.
Therefore,
\Eq{*}{
   H(t_1)\cap H(t_2)\subseteq\bigcup\{H(t)\mid H(t)\subset H(t_1)\cap H(t_2)\}.}
The reversed inclusion is trivial.

Assume that $X,Y$ are topological spaces and $H$ is strongly surjective on $X$. Consider a function $x\colon X\to Y$
fulfilling our requirement. Fix $t_0\in X$ and let $V\subset Y$ be a neighborhood of $x(t_0)$. Since $H$ is strongly
surjective, $t_0\in H(t)^{\circ}$ for some $t\in X$. The restriction $x\rest_{H(t)}$ is continuous, thus there exists
a neighborhood $W\subseteq H(t)^{\circ}$ of $t$ such that $x(W)=x\rest_{H(t)}(W)\subset V$. On the other hand, $W$ can
be represented as $W=U\cap H(t)^{\circ}$, where $U$ is open (in the original topology) and contains $t_0$. Thus
$W\subset X$ is a neighborhood of $t_0$ in the original topology, as well. Therefore, $x\in\CC(X,Y)$.
\end{proof}

Consider the space $\B(X,B)$ of all bounded maps from a nonempty set $X$ to a metric space $(B,d)$. In what follows, we will
equip this space by a family of equivalent norms parametrized by admissible weight functions. A function $p\colon X\to]0,+\infty[$
is termed an \emph{admissible weight function on $X$} if it satisfies $0<\inf_X p\leq\sup_X p<+\infty$; the collection
of such functions is denoted by $\W(X)$. For an arbitrary $p\in\W(X)$ and $x,y\in\B(X,B)$, define
\Eq{*}{
  d_{p}(x,y):=\sup_{t\in X}p(t)d\bigl(x(t),y(t)\bigr).
}
The following result summarizes the properties of the function $d_p$ which will play a key role in the renorming processes.

\Prp{renor}{Let $X$ be a nonempty set and $(B,d)$ be a metric space. Then $\{d_p\mid p\in\W(X)\}$ is a family of pairwise
equivalent metrics on the space $\B(X,B)$. In addition, if $(B,d)$ is complete, then $(\B(X,B),d_p)$ is also complete for
all $p\in\W(X)$. \\ Furthermore, the space $\CC(X,B)$ of all continuous maps from a compact topological space $X$ to a complete
metric space $(B,d)$ is a complete subspace of $(\B(X,B),d_p)$ for all $p\in\W(X)$.}

\begin{proof}
It is elementary to see that $d_p$ is a metric on $\B(X,B)$ for all $p\in\W(X)$. For simplicity, the constant weight function
$p(t)=1$ on $X$ will be denoted by $\mathbf{1}$. One can verify that
\Eq{*}{
   \inf_Xp\cdot d_{\mathbf{1}}\leq d_p\leq \sup_Xp\cdot d_{\mathbf{1}},
}
which proves that $d_p$ is equivalent to $d_{\mathbf{1}}$ for all $p\in\W(X)$. Hence $\{d_p\mid p\in\W(X)\}$ is a family of
pairwise equivalent metrics.

Let $(B,d)$ be a complete metric space. In view of the equivalence of the metrics $d_p$, it is sufficient to show that
$(\B(X,B),d_{\mathbf{1}})$ is a complete metric space. 

As previously, denote the set of all functions from $X$ to $B$ by $B^X$. Then, the Cauchy criterion of uniform convergence
holds in $B^X$: A sequence $(x_n)$ tends to $x\in B^X$ in the supremum distance $d_{\mathbf{1}}$ if and only if, for all
$\eps>0$ there exists $\de>0$ such that
\Eq{*}{
 d\bigl(x_n(t),x_m(t)\bigr)<\eps
}
holds, whenever $n,m>\de$ and $t\in X$. 

Consider now a Cauchy sequence in $(\B(X,B),d_{\mathbf{1}})$. This sequence fulfills the Cauchy criterion and hence converges
uniformly to some element of $B^X$. The triangle inequality guarantees that this element belongs to $\B(X,B)$, which then
yields completeness.

If $X$ is a compact topological space, then $\CC(X,B)$ is a linear subspace of $\B(X,B)$. Therefore any Cauchy
sequence $(x_n)$ of $\CC(X,B)$ is a Cauchy sequence also in $\B(X,B)$. By the previous part, $(x_n)$ tends
to some element $x\in\B(X,B)$ in the supremum distance. Now fix $t_0\in H$ arbitrarily. By the triangle inequality,
\Eq{*}{
 d\bigl(x(t),x(t_0)\bigr)\le 2d_{\mathbf{1}}(x,x_n)+d\bigl(x_n(t),x_n(t_0)\bigr)
}
holds for any $t\in X$. This estimation gives the continuity of $x$ at $t_0$. Therefore $\CC(X,B)$ is complete.
\end{proof}

In a part of the investigations, we will use the Bochner integral \cite{Boc33}. For convenience, we recall its
definition and its most important properties based on Yosida's book \cite{Yos95}. Let $(X,\Sigma,\mu)$ be a measure
space and let $(B,\|\cdot\|)$ be a Banach space. Consider a simple function $x\colon X\to B$ of the form
\Eq{*}{
 x(t):=\sum_{k=1}^n\chi_{E_k}(t)b_k,
}
where $E_1,\dots,E_n$ are pairwise disjoint members of the $\sig$-algebra $\Sigma$, the elements $b_1,\dots,b_n$
belong to $B$, and $\chi_E$ is the characteristic function of $E$. If $\mu(E_k)$ is finite whenever $b_k\neq 0$,
then $x$ is called \emph{Bochner integrable}, and its Bochner integral is defined by
\Eq{*}{
 \int_X xd\mu:=\sum_{k=1}^n\mu(E_k)b_k.
}

A measurable function $x\colon X\to B$ is \emph{Bochner integrable}, if there exists a sequence of integrable simple
functions $(x_n)$ such that
\Eq{*}{
 \lim_{n\to\infty}\int_X\|x-x_n\|d\mu=0,
}
where the integral on the left-hand side is the usual Lebesgue integral. In this case, the Bochner integral of $x$
is given by
\Eq{*}{
 \int_X xd\mu:=\lim_{n\to\infty}\int_X x_nd\mu.
}

It can be shown that the definitions of Bochner integrability and the Bochner integral are independent on the choice
of the approximating sequence. The Bochner integral shares many properties with the Lebesgue integral: It is linear,
$\sig$-additive and fulfills the triangle inequality.

Bochner's majorant condition for integrability plays a distinguished role. A function $f\colon X\to B$ is called
\emph{Bochner-measurable} if it is equal $\mu$-almost everywhere to a function $g$ taking values in a separable
subspace $L$ of $B$, such that the inverse image $g^{-1}(V)$ of every open set $V$ in $B$ belongs to $\Sigma$.
Bochner's criterion states that a Bochner-measurable function $x\colon X\to B$ is Bochner integrable if and only
if
\Eq{*}{
 \int_X\|x\|d\mu <\infty.
}

The last proposition gives a sufficient condition under which continuous functions are Bochner integrable. Its
statement turns out to be crucial in our investigations.

\Prp{Bochner}{If $X$ is a compact topological space with a finite Borel measure and $B$ is a Banach space, then
$\CC(X,B)$ consists of Bochner integrable maps.}

\begin{proof}
Let $x\in\CC(X,B)$ be arbitrary. Since $x$ is continuous and $X$ is compact, $x(X)$ is compact. In particular, $x(X)$ is
completely bounded and thus contains a countable dense subset $D$. The linear hull of $D$ provides a separable subspace $L$:
The rational linear combinations of $D$ is a countable dense subset in $L$. Hence the range of $x$ is contained in a separable
linear subspace of $B$. By the continuity of $x$, the inverse image $x^{-1}(V)$ of any open set $V$ in $B$ is open in $X$.
That is, $x^{-1}(V)$ belongs to the underlying Borel $\sig$-algebra. Therefore $x$ is Borel-measurable.

Using the compactness of $X$, the continuity of $x$, and the finiteness of the Borel measure $\mu$, we arrive at
\Eq{*}{
 \int_X\|x\|d\mu\le\int_X\|x\|_{\infty}d\mu=\|x\|_{\infty}\mu(X)<\infty.
}
Thus the desired statement follows from Bochner's majorant condition.
\end{proof}

Assume that $X$ is a topological space with a Radon measure $\mu$. A relation $H\subseteq X^2$ is called \emph{$\mu$-continuous
at a point $t_0\in X$} if, for all $\eps>0$, there exists a neighborhood $U$ of $t_0$ such that
\Eq{*}{
 \mu\bigl((H(t)\setminus H(t_0))\cup(H(t_0)\setminus H(t))\bigr)<\eps
}
whenever $t\in U$. If $H$ is $\mu$-continuous at each point of $X$, then we say that $H$ is \emph{$\mu$-continuous}. 
%For simplicity, $\mu$-continuity at $t_0$ will be expressed as
%\Eq{*}{
% \lim_{t\to t_0}\mu\bigl((H(t)\setminus H(t_0))\cup(H(t)\setminus %H(t_0))\bigr)=0.
%}

As it is well-known, the integral of an integrable function is continuous at its upper limit. The next proposition
extends this fact and, what is more important, will justify those integral equations which we are going to study.

\Prp{cont}{Let $X$ be a topological space with a Radon measure $\mu$, and let $H\subset X^2$ be a transitive, compact
valued, strongly surjective, and $\mu$-continuous relation. If $B$ is a Banach space and $R\colon H\to B$ is continuous,
then
\Eq{*}{
 \Phi(t):=\int_{H(t)}R(t,s)d\mu(s)
}
defines a continuous map $\Phi\colon X\to\R$.}

\begin{proof}
Note that the definition of $\Phi$ makes sense in view of \prp{Bochner}. Fix $t_0\in X$. The triangle inequality for the
Bochner integral guarantees that
\Eq{*}{\hspace{-5cm}
 \left\|\int_{H(t)}R(t,s)d\mu(s)-\int_{H(t_0)}R(t_0,s)d\mu(s)\right\|
}
\Eq{*}{
 \le&\left\|\int_{H(t)}R(t,s)d\mu(s)-\int_{H(t)\cap H(t_0)}R(t,s)d\mu(s)\right\|\\
   &+\left\|\int_{H(t)\cap H(t_0)}R(t,s)d\mu(s)-\int_{H(t)\cap H(t_0)}R(t_0,s)d\mu(s)\right\|\\
   &+\left\|\int_{H(t)\cap H(t_0)}R(t_0,s)d\mu(s)-\int_{H(t_0)}R(t_0,s)d\mu(s)\right\|\\
 \le&\int_{H(t)\setminus H(t_0)}\left\|R(t,s)\right\|d\mu(s)\\
   &+\int_{H(t)\cap H(t_0)}\left\|R(t,s)-R(t_0,s)\right\|d\mu(s)\\
   &+\int_{H(t_0)\setminus H(t)}\left\|R(t_0,s)\right\|d\mu(s).
}
Here the last term tends to zero as $t\to t_0$ by the absolute continuity of the integral and by the $\mu$-continuity of $H$
at $t_0$. Next we prove the same property of the first term by showing the boundedness of the integrand at a neighborhood of
$t_0$. Since $H$ generates a strongly surjective map, there exists $t^*\in X$ such that $t_0\in H(t^*)^{\circ}$. By transitivity,
$H(t)\subseteq H(t^*)$ if $t\in H(t^*)$.
Thus,
\Eq{*}{
 \sup_{s\in H(t)}\|R(t,s)\|\le\sup\{\|R(t,s)\|\colon t,s\in H(t^*)\}<+\infty,
}
since the right-hand side is the continuous image of a compact set. Using the $\mu$-continuity, we arrive at the desired
statement.

Finally, we show that the middle term tends to zero as $t\to t_0$. Clearly, it is sufficient to prove that
\Eq{*}{
 \lim_{t\to t_0}\sup_{s\in H(t_0)}\|R(t,s)-R(t_0,s)\|=0.
}
Let $\eps>0$ be arbitrary. If $s\in H(t_0)$, then $H(s)\subseteq H(t_0)$. On the other hand, by the continuity of $R$, there
exists a neighborhood $U_s$ of $t_0$ and a neighborhood $V_s$ of $s$, such that
\Eq{*}{
 \|R(t,\sig)-R(t_0,s)\|<\frac{\eps}{2}
}
whenever $(t,\sig)\in (U_s\times V_s)\cap H$. The family of $s$-neighborhoods $\{V_s\mid s\in H(t_0)\}$ is an open cover for
the compact set $H(t_0)$. Thus $H(t_0)\subseteq V_{s_1}\cup\dots\cup V_{s_m}$ holds with suitable $s$-neighborhoods. Define
$U:=U_{s_1}\cap\dots\cap U_{s_m}$. Then, $U$ is a neighborhood of $t_0$. For $(t,\sig)\in (U\times H(t_0))\cap H$, there exists
an index $j\in\{1,\dots,m\}$ such that $\sig\in V_{s_j}$. Hence
\Eq{*}{
 \|R(t,\sig)-R(t_0,\sig)\|\le\|R(t,\sig)-R(t_0,s_j)\|+\|R(t_0,s_j)-R(t_0,\sig)\|<\frac{\eps}{2}+\frac{\eps}{2}=\eps.
}
This completes the proof.
\end{proof}

Observe that the existence of a compact-valued strongly surjective relation $H\subseteq X^2$ has a serious consequence: the
underlying topological space $X$ must be locally compact. Although this fact will not be stated explicitly, our results remain
true in such spaces.

Let $X$ be a topological space with a Radon measure $\mu$ and let $H\subseteq X^2$ be a reflexive, transitive, compact-valued,
strongly surjective, and $\mu$-continuous relation. Then, by \prp{cont}, the map $\Lambda_{H,\mu}$ defined by
\Eq{*}{
  \big(\Lambda_{H,\mu}x\big)(t):= \int_{H(t)}xd\mu
}
is a linear selfmap of the space $\CC(X,\R)$. This we will be called the \emph{core map associated to the pair $(H,\mu)$}. In this
context, it is also natural to introduce the \emph{spectral radius function of $\Lambda_{H,\mu}$} by
\Eq{*}{
  \rho_{H,\mu}(t):=&\limsup_{k\to\infty}\big(\big(\Lambda_{H,\mu}^k \mathbf{1}\big)(t)\big)^{\frac1k}\\
  =&\limsup_{k\to\infty}\Bigg(\int_{H(t)}\bigg(\int_{H(s_1)}\dots\bigg(\int_{H(s_{k-1})} d\mu(s_k) \bigg)\dots d\mu(s_2)\bigg)
   d\mu(s_1)\Bigg)^{\frac1k}.
}
The spectral radius function is monotonic in the following sense: If $s\in H(t)$, then the transitivity of $H$ implies that
$H(s)\subseteq H(t)$, and hence, $\rho_{H,\mu}(s)\leq \rho_{H,\mu}(t)$.

\Prp{CH}{Let $X$ be a topological space with a Radon measure $\mu$, and let $H\subset X^2$ be a reflexive, transitive, compact
valued, strongly surjective, and $\mu$-continuous relation. Let $L_0\geq0$ and $t_0\in X$ such that
\Eq{*}{
  L_0\cdot\rho_{H,\mu}(t_0)<1.
}
Then, the integral equation 
\Eq{IEl}{
  \ell(t)=1+L_0\int_{H(t)}\ell(s) d\mu(s)
}
has a positive solution $\ell$ in $\colon H(t_0)\to\R_+$.}

\begin{proof}
Define the sequence of real valued functions $(\ell_n)_{n=0}^\infty$ on $H(t_0)$ by
\Eq{*}{
 \ell_{n}(t):=\sum_{k=0}^nL_0^k\cdot\big(\Lambda_{H,\mu}^k \mathbf{1}\big)(t).
}
Then, $(\ell_n)$ is a nondecreasing sequence whose members belong to the space $\CC(H(t_0),\R)$ by \prp{cont}. Moreover $\ell_0=\mathbf 1$,
which implies that $\mathbf{1}=\ell_0\le \ell_n$. It is also easy to see that, for all $t\in H(t_0)$ and $n\in\N$,
\Eq{BPI}{
 \ell_{n}(t)
 =1+L_0\cdot\big(\Lambda_{H,\mu}\ell_{n-1}\big)(t)
 =1+L_0\int_{H(t)}\ell_{n-1}(s)d\mu(s),
}
which shows that $(\ell_n)$ is a Banach--Piccard iteration sequence.

The Cauchy--Hadamard Theorem and the assumption $L_0\cdot\rho_{H,\mu}(t_0)<1$ guarantee that the series
\Eq{*}{
 \sum_{k=0}^{\infty}L_0^k\cdot\big(\Lambda_{H,\mu}^k \mathbf{1}\big)(t_0)
}
is convergent. On the other hand, for $t\in H(t_0)$, the transitivity of $H$ implies that $H(t)\subseteq H(t_0)$, hence, for all $k\in\N$,
\Eq{*}{
\big(\Lambda_{H,\mu}^k\mathbf{1}\big)(t)
=\int_{H(t)}\big(\Lambda_{H,\mu}^{k-1}\mathbf{1}\big)(s)d\mu(s)
\le\int_{H(t_0)}\big(\Lambda_{H,\mu}^{k-1}\mathbf{1}\big)(s)d\mu(s)
=\big(\Lambda_{H,\mu}^k\mathbf{1}\big)(t_0). 
}
Therefore, the Weierstrass convergence theorem yields that
\Eq{*}{
 \sum_{k=0}^{\infty}L_0^k\cdot\big(\Lambda_{H,\mu}^k \mathbf{1}\big)(t)
}
is uniformly convergent for $t\in H(t_0)$. The members of $(\ell_n)$ are continuous, therefore the pointwise limit function
$\ell:=\lim_{n\to\infty}\ell_n$ is also continuous on $H(t_0)$. The inequality $\mathbf{1}\leq \ell_n$ implies that $\ell$
is positive everywhere. Finally, upon taking the limit $n\to\infty$ in \eqref{BPI}, it follows that $\ell$ satisfies the
integral equation of the theorem. 
\end{proof}

If the spectral radius function of the core map associated to $(H,\mu)$ is equal to zero at some $t_0\in X$, then as an immediate
consequence of the \prp{CH}, we get that the integral equation \eqref{IEl} has a positive continuous solution on $H(t_0)$. 

\section{Nonlinear integral equations}

Our main results are presented in three theorems. The first one concludes the unique resolvability of nonlinear integral equations
provided that there exists a solution of the corresponding linear homogeneous integral inequality. This assumption makes possible
to apply the renorming technique of Bielecki. Moreover, it can easily be checked in practice via numerical methods.

\Thm{FredVolt}{Let $X$ be a topological space with a Radon measure $\mu$, and let $H\subset X^2$ be a reflexive, transitive,
compact valued, strongly surjective, and $\mu$-continuous relation. Let $B$ be a Banach space, and assume that the continuous
kernel $K\colon H\times B\to B$ fulfills the Lipschitz condition
\Eq{*}{
 \|K\bigl(t,s,x\bigr)-K\bigl(t,s,y\bigr)\|\le L(t,s)\|x-y\|
}
for all $(t,s)\in H$ and $x,y\in B$ with a continuous function $L\colon H\to\R_+$. If $f\in\CC(X,B)$ and, for all $t_0\in X$, the
linear homogeneous integral inequality
\Eq{II}{
 \int_{H(t)}L(t,s)\ell(s)d\mu(s)<\ell(t)
}
has a positive solution $\ell$ in $\colon H(t_0)\to\R$, then the nonlinear integral equation
\Eq{FV}{
 x(t)=f(t)+\int_{H(t)}K\bigl(t,s,x(s)\bigr)d\mu(s)
}
has a unique solution $x$ in $\CC(X,B)$.}

\begin{proof}
Note that the inequality \eqref{II} is correctly formulated by the transitivity of $H$, and that \eqref{FV} makes sense in view of
\prp{Bochner}. Now consider the map $T$ defined by
\Eq{*}{
 (Tx)(t):=f(t)+\int_{H(t)}K(t,s,x(s))d\mu(s).
}
By \prp{cont}, the right-hand side above is a continuous function of $t$. Thus, $T$ is a self-map of the space $\CC(X,B)$.

Fix now $t_0\in X$. If $x,y\in\CC(X,B)$ fulfill $x\rest_{H(t_0)}=y\rest_{H(t_0)}$, then by the transitivity of $H$, for all $t\in H(t_0)$,
we obtain
\Eq{*}{
 (Tx)\rest_{H(t_0)}(t)
  &=f(t)+\int_{H(t)}K(t,s,x(s))d\mu(s)\\
   &=f(t)+\int_{H(t)}K(t,s,x\rest_{H(t_0)}(s))d\mu(s)\\
    &=f(t)+\int_{H(t)}K(t,s,y\rest_{H(t_0)}(s))d\mu(s)\\
     &=f(t)+\int_{H(t)}K(t,s,y(s))d\mu(s)=(Ty)\rest_{H(t_0)}(t).
}
This shows that $T$ is restrictable to $\CC(H(t_0),B)$. Next we prove that this restriction, denoted by $T$ as well, has a unique fixed
point in $\CC(H(t_0),B)$. Let $\ell\colon H(t_0)\to\R_+$ be a positive and continuous solution of \eqref{II}. By the compactness of $H(t_0)$
and by \prp{cont} again,
\Eq{*}{
 q:=\max_{t\in H(t_0)}\frac{1}{\ell(t)}\int_{H(t)}L(t,s)\ell(s)d\mu(s)<1.
}
For $x,y\in \CC(H(t_0),B)$, we have $Tx,Ty\in \CC(H(t_0),B)$. Applying the Lipschitz-condition, with the notation $p:=1/\ell$, we get that,
for all $t\in H(t_0)$,
\Eq{*}{
 p(t)\|(Tx)(t)-(Ty)(t)\| 
  & \le p(t) \int_{H(t)}\|K(t,s,x(s))-K(t,s,y(s))\|d\mu(s)\\
  & \le p(t) \int_{H(t)}L(t,s)\|x(s)-y(s)\|d\mu(s)\\
  & \le p(t) \int_{H(t)}\frac{L(t,s)}{p(s)}\|x(s)-y(s)\|p(s)d\mu(s)\\
  & \le \bigg(p(t)\int_{H(t)}\frac{L(t,s)}{p(s)}d\mu(s)\bigg)\cdot \|x-y\|_p\\
  & = \bigg(\frac{1}{\ell(t)}\int_{H(t)}L(t,s)\ell(s)d\mu(s)\bigg)\cdot \|x-y\|_p\\
  & \le q\|x-y\|_p.
}
Taking supremum in $t\in H(t_0)$ in the initial term, we arrive at $\|Tx-Ty\|_p\le q\|x-y\|_p$. This means that the restriction of $T$
to $\CC(H(t_0),B)$ is a contraction in the $p$-norm. Thus $T$ has a unique fixed point in $\CC(H(t_0),B)$ by the Contraction
Principle. Finally, \prp{rel} and \prp{solex} complete the proof.
\end{proof}

It is important to observe that, under these assumptions, our theorem implies that the integral equation
\Eq{IE}{
 \ell(t)=1+\int_{H(t)}L(t,s)\ell(s)d\mu(s)
}
has a solution $x\in\CC(X,\R)$. On the other hand, consider the sequence $(\ell_n)$ determined by the Banach--Piccard iteration
\Eq{*}{
  \ell_0(t)=1,\qquad \ell_{n}(t):=1+\int_{H(t)}L(t,s)\ell_{n-1}(s)d\mu(s).
}
Then, $(\ell_n)$ is nondecreasing with respect to the pontwise ordering, it converges to $\ell$, and the convergence is uniform
on $H(t_0)$ for all $t_0\in X$. Therefore, $1=\ell_0\leq\ell$ shows that $\ell$ is a \emph{positive continuous} solution of \eqref{IE}
and thus also of \eqref{II} over the \emph{entire} set $X$.

The next result is a global existence and uniqueness theorem for the solvability of nonlinear integral equations. The role of the
inequality \eqref{II} is hidden: Instead, we use the spectral radius function.

\Thm{FredVolt+}{Let $X$ be a topological space with a Radon measure $\mu$, and let $H\subset X^2$ be a reflexive, transitive, compact
valued, strongly surjective, and $\mu$-continuous relation. Let $B$ be a Banach space, and assume that the continuous kernel
$K\colon H\times B\to B$ fulfills the Lipschitz condition
\Eq{*}{
 \|K\bigl(t,s,x\bigr)-K\bigl(t,s,y\bigr)\|\le L(t)\|x-y\|
}
for all $(t,s)\in H$ and $x,y\in B$ with a continuous function $L\colon X\to\R_+$. If $f\in\CC(X,B)$ and the spectral radius function
$\rho_{H,\mu}$ is identically zero on $X$, then the nonlinear integral equation \eqref{FV} has a unique solution $x$ in $\CC(X,B)$.}

\begin{proof} 
Obviously, the Lipschitz condition of this theorem implies the weaker Lipschitz condition of \thm{FredVolt}. In order to draw the conclusion
of this theorem, it is enough to verify the existence of a positive continuous solution of the integral inequality \eqref{II} for all $t_0\in X$.
For this goal, it is sufficient to prove the solvability of the integral inequality
\Eq{*}{
  L_0\int_{H(t)}\ell(s)d\mu(s)<\ell(t),
}
where $L_0=\sup_{t\in H(t_0)}L(t)$. This assertion, however, directly follows from $\rho_{H,\mu}(t_0)=0$ and \prp{CH} because any solution
of \eqref{IEl} is also a solution of the above inequality.
\end{proof}

Finally, we state a theorem which allows multivariable kernels. Let us emphasize that, due to this property, it allows even retardations in
the nonlinear integral equation.

\Thm{MultFredVolt}{Let $X$ be a topological space with a Radon measure $\mu$, and let $H\subset X^2$ be a reflexive, transitive,
compact valued, strongly surjective, and $\mu$-continuous relation. Let $B$ be a Banach space, and assume that the continuous
kernel $K\colon H\times B^n\to B$ fulfills the Lipschitz condition
\Eq{*}{
 \|K\bigl(t,s,x_1,\dots,x_n\bigr)-K\bigl(t,s,y_1,\dots,y_n\bigr)\|\le\sum_{k=1}^nL_k(t,s)\|x_k-y_k\|
}
for all $(t,s)\in H$ and $x_k,y_k\in B$ with continuous functions $L_k\colon H\to\R_+$. If $f\in\CC(X,B)$, the functions
$\phi_1,\dots,\phi_n\in\CC(X,X)$ satisfy $\phi_k\circ H\subseteq H$ for all $k\in\{1,\dots,n\}$, and, for all $t_0\in X$, the linear
homogeneous integral inequality
\Eq{*}{
 \sum_{k=1}^n\int_{H(t)}L_k(t,s)\ell(\phi_k(s))d\mu(s)<\ell(t)
}
has a positive solution $\ell$ in $\colon H(t_0)\to\R$, then the nonlinear retarded integral equation
\Eq{*}{
 x(t)=f(t)+\int_{H(t)}K\bigl(t,s,x(\phi_1(s)),\dots,x(\phi_n(s))\bigr)d\mu(s)
}
has a unique solution $x$ in $\CC(X,B)$.}

Obviously, \thm{MultFredVolt} implies \thm{FredVolt}. However, the proof of the above result is completely similar to that of \thm{FredVolt},
therefore it is omitted.

\section{Applications to Fredholm-type equations}

If $X$ is a compact topological space with a Radon measure $\mu$, then $H=X^2$ is a reflexive, transitive, and $\mu$-continuous
relation on $X$. Using this easy observation, \thm{FredVolt} reduces to the next Fredholm-type
result:

\Thm{Fredholm}{Let $X$ be a compact topological space with a Radon measure $\mu$, and let $B$ be a Banach space. Assume that the
continuous kernel $K\colon X^2\times B\to B$ fulfills the Lipschitz condition
\Eq{*}{
  \|K(t,s,x)-K(t,s,y)\|\leq L(t,s)\|x-y\|
}
for all $t,s\in X$ and $x,y\in B$ with a continuous function $L\colon X^2\to\R_+$. If the linear homogeneous integral inequality
\Eq{*}{
 \int_XL(t,s)\ell(s)d\mu(s)<\ell(t)
}
has a positive solution $\ell$ in $\colon X\to\R_+$, then the nonlinear Fredholm-type equation
\Eq{Fred}{
 x(t)=f(t)+\int_XK(t,s,x(s))d\mu(s)
}
has a unique solution $x$ in $\CC(X,B)$.}

Not claiming completeness, we sketch two consequences of this result. The first corollary is a special case of \thm{Fredholm} if
$\ell\equiv 1$. In the second one, we assume that the Lipschitz modulus has a product form.

\Cor{F1}{Let $X$ be a compact topological space with a Radon measure $\mu$, and let $B$ be a Banach space. Assume that the continuous
functions $K\colon X^2\times B\to B$ and $L\colon X^2\to\R_+$ fulfill the conditions
\Eq{*}{
  \|K(t,s,x)-K(t,s,y)\|\leq L(t,s)\|x-y\|\quad\mbox{ and }\quad
  \int_XL(t,s)d\mu(s)<1
}
for all $t,s\in X$ and $x,y\in B$. Then the nonlinear integral equation \eqref{Fred} has a unique continuous solution $x\colon X\to B$.}

\Cor{F2}{Let $X$ be a compact topological space with a Radon measure $\mu$, and let $B$ be a Banach space. Assume that the continuous
functions $K\colon X^2\times B\to B$ and $L_1,L_2\colon X\to\R_+$ fulfill the conditions
\Eq{*}{
  \|K(t,s,x)-K(t,s,y)\|\leq L_1(t)L_2(s)\|x-y\|\quad\mbox{ and }\quad
  \int_XL_1(s)L_2(s)d\mu(s)<1
}
for all $t,s\in X$ and $x,y\in B$. Then the nonlinear integral equation \eqref{Fred} has a unique continuous solution $x\colon X\to B$.}

\begin{proof}
First choose $c>0$ so that 
\Eq{*}{
 \int_H(L_1(s)+c)L_2(s)d\mu(s)<1
}
be valid. Then, all conditions of \thm{Fredholm} are satisfied with $L(t,s):=(L_1(t)+c)L_2(s)$ and $\ell:=L_1+c$. Indeed,
\Eq{*}{
  \int_H L(t,s)\ell(s)d\mu(s)
  &= \int_H (L_1(t)+c)L_2(s)(L_1(s)+c)d\mu(s)\\
  &= (L_1(t)+c)\int_H L_2(s)(L_1(s)+c)d\mu(s)\\
  &<L_1(t)+c=\ell(t).
}
Thus the statement follows from \thm{Fredholm}.
\end{proof}

The particular cases of \cor{F1} and \cor{F2}, when the kernel is the linear transform of the unknown
function, may also be mentioned:

\Cor{F3}{Let $X$ be a compact topological space with a Radon measure $\mu$. If $f\colon X\to\R^m$ and $A\colon X^2\to\R^{m\times m}$
are continuous and satisfy
\Eq{*}{
 \int_X \|A(t,s)\| d\mu(s)<1,
}
then the inhomogeneous linear integral equation
\Eq{*}{
 x(t)=f(t)+\int_XA(t,s)x(s)d\mu(s) 
}
has a unique continuous solution $x\colon X\to\R^m$.}

\Cor{F4}{Let $X$ be a compact topological space with a Radon measure $\mu$. If $f\colon X\to\R^m$, further $A_1\colon X\to\R^{m\times k}$
and $A_2\colon X\to\R^{k\times m}$ are continuous such that 
\Eq{*}{
 \int_X \|A_1(s)\|\|A_2(s)\| d\mu(s)<1,
}
then the inhomogeneous linear integral equation
\Eq{*}{
 x(t)=f(t)+\int_XA_1(t)A_2(s)x(s)d\mu(s)
}
has a unique continuous solution $x\colon X\to\R^m$.}

Observe that the original result of Fredholm follows from \cor{F3} in the special setting when $X$ is a compact interval, $\mu$ is the
Lebesgue measure, and the kernel takes real values.

\section{Applications to Volterra-type equations}

The standard exposition of Volterra's result proves that a suitable iterate of the map defined via Volterra's equation is a contraction.
Instead of the standard approach, we use the renorming technique and obtain a more general result. Now the solvability of the corresponding
homogeneous linear inequality \eqref{II} remains hidden.

Consider the standard partial order $\leq$ on $\R^n$ induced by the nonnegative orthant $[0,\infty[^n$: For $a,b\in\R^n$, the inequality
$a\leq b$ means that the coordinates of $b-a$ are nonnegative. In this case, we define the $n$-dimensional interval (rectangle)
$[a,b]\subseteq\R^n$ by
\Eq{*}{
  [a,b]:=\{u\in\R^n\mid a\leq u\leq b\}.
}
We say that a set $D\subseteq[0,\infty[^n$ is \emph{rectangular (with respect to the origin)} if, for all $u\in D$, the rectangle $[0,u]$
is contained in $D$ and the set $[0,\infty[^n\setminus D$ is closed in $\R^n$. 

\Thm{Volterra}{Let $D\subseteq [0,\infty[^n$ be a rectangular set, let $\Delta(D):=\{(t,s)\mid t\in D,\,s\in[0,t]\}$, and let $B$ be a
Banach space. Assume that the continuous kernel $K\colon\Delta(D)\times B\to B$ fulfills the Lipschitz condition
\Eq{*}{
 \|K\bigl(t,s,x\bigr)-K\bigl(t,s,y\bigr)\|\le L(t)\|x-y\|
}
for all $(t,s)\in\Delta(D)$ and $x,y\in B$ with a continuous function $L\colon D\to\R_+$. If $f\colon D\to B$ is a continuous function,
then the nonlinear Volterra-type equation
\Eq{Volt}{
 x(t)=f(t)+\int_{[0,t]} K\bigl(t,s,x(s)\bigr)ds
}
has a unique solution $x$ in $\CC(D,B)$.}

\begin{proof}
Let $X=D$ be equipped with the Euclidean subspace topology and set $H:=\Delta(D)$. Then $H(t)=[0,t]$ for all $t\in D$, showing that the
values of $H$ are compact. Furthermore, $H$ is reflexive, transitive, strongly surjective, and continuous with respect to the $n$-dimensional
Lebesgue measure $\lambda$. 

In what follows, we show that the spectral radius function $\rho_{H,\lambda}$ is identically zero on $D$. To accomplish this goal, first
we prove by induction on $k$ that, for all $t\in D$,
\Eq{k}{
  (\Lambda_{H,\lambda}^k \mathbf{1})(t)=\frac{\lambda([0,t])^k}{(k!)^n}.
}
For $k=1$, we have
\Eq{*}{
  (\Lambda_{H,\lambda}\mathbf{1})(t)=\int_{[0,t]}1d\lambda(s)=\lambda([0,t]).
}
Now assume that \eqref{k} holds for some $k$, and let $t=(t_1,\dots,t_n)$. Then, by Fubini's theorem,
\Eq{*}{
  (\Lambda_{H,\lambda}^{k+1} \mathbf{1})(t)
  &=\int_{[0,t]}(\Lambda_{H,\lambda}^k \mathbf{1})(s)d\lambda(s)
  =\int_{[0,t]}\frac{\lambda([0,s])^k}{(k!)^n}d\lambda(s)\\
  &=\frac{1}{(k!)^n}\int_{[0,t]}(s_1\cdots s_n)^kd\lambda(s_1,\dots,s_n)
  =\frac{1}{(k!)^n}\prod_{i=1}^n\int_0^{t_i} s_i^k ds_i\\
  &=\frac{1}{(k!)^n}\prod_{i=1}^n\frac{t_i^{k+1}}{k+1}
  =\frac{\lambda([0,t])^{k+1}}{((k+1)!)^n}.
}
Finally, using the just have proved equality \eqref{k}, for the spectral radius function $\rho_{H,\lambda}$, we obtain
\Eq{*}{
  \rho_{H,\lambda}(t)
  =\limsup_{k\to\infty}
   \big((\Lambda_{H,\lambda}^k \mathbf{1})(t)\big)^{\frac1k}
  =\limsup_{k\to\infty}
   \frac{\lambda([0,t])}{\big(\sqrt[k]{k!}\big)^n}=0.
}
Thus, the assertion directly follows from \thm{FredVolt+}. 
\end{proof}

As the most important applications to \thm{Volterra}, we present here two Corollaries. The first one extends the classical
result of Volterra. In this extension, the Lipschitz property of the kernel can be checked directly, therefore we omit the
details of the proof.

\Cor{V1}{Let $D\subseteq [0,\infty[^n$ be a rectangular set and let $\Delta(D):=\{(t,s)\mid t\in D,\,s\in[0,t]\}$. If $f\colon D\to\R^m$
and $A\colon\De(D)\to\R^{m\times m}$ are continuous functions, then the inhomogeneous linear Volterra equation
\Eq{*}{
 x(t)=f(t)+\int_{[0,t]}A(t,s)x(s)ds
}
has a unique continuous solution $x\colon D\to\R^m$.}

The second application of \thm{Volterra} is an existence and uniqueness result for a boundary value problem for a special kind
of partial differential equation. In the particular case when $n=1$, it reduces to the Global Existence and Uniqueness Theorem,
while for $n=2$, it extends the one-dimensional Wave Equation. 

\Cor{V2}{Let $D\subseteq [0,\infty[^n$ be a rectangular set and denote 
\Eq{*}{
 D_0:=\{(t_1,\dots,t_n)\in D\mid t_1\cdots t_n=0\}.
}
Let $F\colon D\times\R^m\to\R^m$ be a continuous map which fulfills the Lipschitz condition
\Eq{*}{
 \|F(t,x)-F(t,y)\|\le L(t)\|x-y\|
}
for all $t\in D$ and $x,y\in\R^m$ with a continuous function $L\colon D\to\R_+$ and let $\phi\colon D_0\to\R^m$ be a continuous function
such that partial derivatives $\partial_n\phi$, $\partial_{n-1}\partial_n\phi$, \dots, $\partial_2\cdots\partial_n\phi$ exist on $D_0$.
Then there exists precisely one continuous solution $x\colon D\to\R^m$ of the Cauchy problem
\Eq{CP}{
 \partial_1\cdots\partial_n x(t)=F\bigl(t,x(t)\bigr)\quad\mbox{on $D$},
  \qquad x(t)=\phi(t)\quad\mbox{on $D_0$}.
}}

\begin{proof}
Denote the set of singular and diagonal $n\times n$ matrices with entries in $\{0,1\}$ by $\Pi_n$.  First we show that \eqref{CP} is
equivalent to the integral equation
\Eq{IEV}{
 x(t)=f(t)+\int_{[0,t]} F\bigl(s,x(s)\bigr)ds,
}
where $f\colon D\to\R^m$ is defined by
\Eq{*}{
 f(t)=\sum_{P\in\Pi_n}(-1)^{n-1-\rank(P)}\phi(Pt).
}

For a fixed $\tau\in[0,\infty[$ and $k\in\{1,\dots,n\}$ such that $D\cap(D-\tau e_k)\neq\emptyset$, introduce the
difference operator $\Delta_{k;\tau}\colon\CC(D,\R^m)\to\CC(D\cap(D-\tau e_k),\R^m)$ by
\Eq{*}{
  (\Delta_{k;\tau} x)(t)=x(t+\tau e_k)-x(t),
}
where $e_k$ stands for the $k$th member of the standard base of $\R^n$.

Assume that $x\colon D\to\R^m$ is a continuous solution of the Cauchy problem \eqref{CP}. Then $x$ is partially
differentiable with respect to its $n$th variable on $D$. Similarly, $\partial_n x$ is partially differentiable
with respect to its $(n-1)$st variable on $D$. Finally, $\partial_2\cdots\partial_n x$ is  partially differentiable
with respect to its first variable on $D$ and \eqref{CP} hold on the indicated domains. The first equality shows that
$\partial_1\cdots\partial_n x$ is continuous on $D$. 

Let $t=(t_1,\dots,t_n)\in D$ be fixed. Integrating the first equality in \eqref{CP} side by side with respect to the
first variable on the interval $[0,t_1]$ and using the Newton--Leibniz formula in the first variable, we obtain
\Eq{*}{
  (\Delta_{1;t_1}\partial_2\cdots\partial_n x)(0,t_2,\dots,t_n)
  =\int_0^{t_1} F(s_1,t_2,\dots,t_n,x(s_1,t_2,\dots,t_n))ds_1.
}
Applying the same process on the forthcoming variables, finally we arrive at
\Eq{*}{
  (\Delta_{1;t_1}\cdots\Delta_{n;t_n} x)(0,\dots,0) &=
   \int_0^{t_n}\cdots\int_0^{t_1} F(s_1,\dots,s_n,x(s_1,\dots,s_n))ds_1\cdots ds_n\\ &=
    \int_{[0,t]}F(s,x(s))d\lambda(s).
}
On the other hand, using the boundary condition of \eqref{CP}, one can see that 
\Eq{*}{
  (\Delta_{1;t_1}\cdots\Delta_{n;t_n} x)(0,\dots,0)=
   x(t)+\sum_{P\in\Pi_n}(-1)^{n-\rank(P)}\phi(Pt)=
    x(t)-f(t).
}

For the converse statement, observe that the partial derivatives $\partial_n f$, $\partial_{n-1}\partial_n f$,
\dots, $\partial_2\cdots\partial_n f$ exist on $D$ by the similar properties of $\phi$. Furthermore, each term in the
definition of $f$ is independent of one of the variables, thus $\partial_1\cdots\partial_n f=0$ on $D$. On the other
hand, by Fubini Theorem,
\Eq{*}{
  x(t)=f(t)+\int_0^{t_n}\cdots\int_0^{t_1} F(s_1,\dots,s_n,x(s_1,\dots,s_n))ds_1\cdots ds_n
}
holds for any continuous solution $x\colon D\to\R^m$ of \eqref{IEV} and for all $t=(t_1,\dots,t_n)\in D$. Differentiating
both sides with respect to $t_n,\dots,t_1$, we obtain the first equality in \eqref{CP}. For the second equality, let
$t=(t_1,\dots,t_n)\in D_0$ be fixed. Then there exists $k\in\{1,\dots,n\}$ such that $t_k=0$. In this case, the above
equation yields $x(t)=f(t)$. Thus it remains to show that $f(t)=\phi(t)$ holds.

Let $E_k\in\Pi_n$ be that matrix whose entries are zero, except for the $k$th member of its diagonal. Let $P_k:=E-E_k$,
where $E$ is the $n\times n$ unit matrix. Consider the transformation $T\colon\Pi_n\to\Pi_n$ given by $T(P):=P\stackrel{.}{+} E_k$,
where $\stackrel{.}{+}$ stands for the addition in the set $\{0,1\}$ modulo 2. Clearly, $T$ is a bijection on $\Pi_n\setminus\{P_k\}$,
and $t_k=0$ ensures $T(P)t=Pt$. Furthermore, the parity of the ranks of $P$ and $T(P)$ are the opposite, and $P_kt=t$. Thus,
\Eq{*}{
  f(t) &=\sum_{P\in\Pi_n}(-1)^{n-1-\rank(P)}\phi(Pt)\\ &=
   \phi(P_kt)+\frac{1}{2}\sum_{P\in\Pi_n\setminus\{P_k\}}\Big((-1)^{n-1-\rank(P)}\phi(Pt)+(-1)^{n-1-\rank(T(P))}\phi(T(P)t)\Big)\\ &=
    \phi(P_kt)+\frac{1}{2}\sum_{P\in\Pi_n\setminus\{P_k\}}\Big((-1)^{n-1-\rank(P)}+(-1)^{n-1-\rank(T(P))}\Big)\phi(Pt)\\ &=
     \phi(P_kt)=\phi(t).
}
Applying \thm{Volterra} to the equivalent integral form of our Cauchy problem, we get the statement of the Corollary. 
\end{proof}

As we have pointed out, the Global Existence and Uniqueness Theorem of ODEs is a direct consequence of \cor{V2}. However, it is also
worth mentioning that the original approach of Bielecki to this result manifests in the proof of \thm{FredVolt}. Indeed, let $H\subseteq\R^2$
be defined by $H(t):=[\tau,t]$ for $\tau<t$, and let $f\equiv\xi\in\R^n$ be fixed initial value. Assume that the Lipschitz-modulus
depends only on the second variable of the kernel. If
\Eq{*}{
 \ell(t):=\exp\left(\int_{\tau}^tL(s)ds\right),
}
then
\Eq{*}{
 \int_{\tau}^tL(s)\ell(s)ds &=
 \int_{\tau}^tL(s)\exp\left(\int_{\tau}^sL(u)du\right)ds=
 \int_{\tau}^t\frac{d}{ds}\exp\left(\int_{\tau}^sL(u)du\right)ds\\ &=
 \exp\left(\int_{\tau}^tL(u)du\right)-1=\ell(t)-1<\ell(t).
}
Thus $\ell$ solves the corresponding inequality of \eqref{II} in this particular setting. Observe also, that the formula \eqref{k} in
\thm{FredVolt+} shows a tight analogue with the standard approach to Volterra's result. However, this standard approach proves that a
suitable iterate of the map defined by the original integral equation is a contraction in the \emph{original} supremum norm.

\section{Applications to Presi\'c-type equations}

Finally we investigate a class of single valued functional equations introduced by Presi\'c in \cite{Pre63a} and
\cite{Pre63b}. An excellent English exposition of his results can be found in \cite{KucChoGer90}. The interaction
of algebraic and analytic aspects of this topic is elaborated in \cite{BesKezHor12} and \cite{BesKez13}. The pure
algebraic feature of the linear Presi\'c-type equation is completely described in \cite{BesKonSza17}. The next result
is a sufficient condition for the unique solvability of the Presi\'c-type equation, with no algebraic restrictions
on the substituting functions. Let us emphasize, that both its philosophy and idea of proof are the same as in
\thm{Fredholm}.

\Thm{Presic}{Let $H$ be a reflexive, compact valued, strongly surjective relation on a topological space $X$, and
let $(B,d)$ be a complete metric space. Assume that $\phi_1,\dots,\phi_n\in\CC(X,X)$ are such that $\phi_k\circ H\subseteq H$
for all $k\in\{1,\dots,n\}$, and $F\in\CC(X\times B^n,B)$ satisfies the Lipschitz condition
\Eq{*}{
 d\bigl(F(t,x_1,\dots,x_n),F(t,y_1,\dots,y_n)\bigr)\le
  \sum_{k=1}^nL_k(t)d(x_k,y_k)
}
for all $t\in X$ and $x_k,y_k\in B$ with continuous functions $L_k\colon X\to]0,+\infty[$. If there exists a
positive function $\ell$ in $\CC(X,\R)$ satisfying
\Eq{*}{
 \sum_{k=1}^nL_k(t)\ell\bigl(\phi_k(t)\bigr)<\ell(t)
}
for all $t\in X$, then there exists a unique solution $f$ in $\CC(X,B)$ of the Presi\'c-type functional
equation
\Eq{Pres}{
 f(t)=F\bigl(t,f(\phi_1(t)),\dots,f(\phi_n(t))\bigl).
}}

\begin{proof}
For $f\in\CC(X,B)$, consider the map $T$ defined by
\Eq{*}{
  (Tf)(t):=F\big(t,f(\phi_1(t)),\dots,f(\phi_n(t))\big).
}
Clearly, $T\colon\CC(X,B)\to\CC(X,B)$. Fix $t\in X$ and assume that two continuous functions $f,g$ fulfill $f\rest_{H(t)}=g\rest_{H(t)}$. 
The inclusion $\phi_k\bigl(H(t)\bigr)\subseteq H(t)$ then implies  
\Eq{*}{
 f\rest_{H(t)}(\phi_k(s))=g\rest_{H(t)}(\phi_k(s))
}
for all $s\in H(t)$ and $k=1,\dots,n$. Therefore,
\Eq{*}{
 Tf\rest_{H(t)}(s)
  &=F\big(s,f\rest_{H(t)}(\phi_1(s)),\dots,f\rest_{H(t)}(\phi_n(s))\big)\\
   &=F\big(s,g\rest_{H(t)}(\phi_1(s)),\dots,g\rest_{H(t)}(\phi_n(s))\big)=Tg\rest_{H(t)}(s).
}
Thus $T$ is restrictable to $\CC(H(t),B)$. Now we prove that every natural restriction, denoted by $T$ as well, has a
unique fixed point. The compactness of $H(t)$, the continuity of $L_k$ and $\ell$, furthermore our assumption imply
\Eq{*}{
  q:=\max_{s\in H(t)}\frac{1}{\ell(s)}\sum_{k=1}^nL_k(s)\ell\bigl(\phi_k(s)\bigr)<1.
}
We claim that $T$ is a $q$-contraction in the complete metric space $(\CC(H(t),B),d_p)$, where the weight function is
given by $p:=1/\ell$. Let $f,g\in\CC(H(t),B)$ be arbitrary. Applying the Lipschitz-condition and the definitions of $q$
and $p$,
\Eq{*}{
  p(s)&d\bigl((Tf)(s),(Tg)(s)\bigr)\\
  &=\frac{1}{\ell(s)}d\bigl(F\big(s,f(\phi_1(s)),\dots,f(\phi_n(s))\bigr),
      F\big(s,g(\phi_1(s)),\dots,g(\phi_n(s))\big)\bigr)\\
  &\leq\frac{1}{\ell(s)}\sum_{k=1}^nL_k(s)d\bigl(f(\phi_k(s)),g(\phi_k(s))\bigr)\\
  &=\frac{1}{\ell(s)}\sum_{k=1}^nL_k(s)\ell(\phi_k(s))\frac{d\bigl(f(\phi_k(s)),g(\phi_k(s)\bigr)}{\ell(\phi_k(s)\bigr)}\\
  &\le\frac{1}{\ell(s)}\sum_{k=1}^nL_k(s)\ell(\phi_k(s))d_p(f,g)\\
  &\leq qd_p(f,g).
}
Taking supremum for $s\in H(t)$ in the left-hand-side term, we arrive at the desired contractivity property. Therefore, $T$ has a
unique fixed point in $\CC(H(t),B)$ by the Contraction Principle. Finally, \prp{rel} and \prp{solex} show that $T$ has a unique fixed
point in $\CC(X,B)$. This fixed point is the unique continuous solution of the Presi\'c-equation.
\end{proof}

The following result is an immediate consequence of the previous theorem.

\Cor{pres1}{Let $H$ be a reflexive, compact valued, strongly surjective relation on a topological space $X$, and
let $B$ be a Banach space. Assume that $\phi_1,\dots,\phi_n\in\CC(X,X)$ fulfill $\phi_k(H(t))\subseteq H(t)$
for all $k\in\{1,\dots,n\}$ and $t\in X$, and $F\in\CC(X\times B^n,B)$ satisfies the Lipschitz condition
\Eq{*}{
 \|F(t,x_1,\dots,x_n)-F(t,y_1,\dots,y_n)\|\le
  \sum_{k=1}^nL_k(t)\|x_k-y_k\|\quad\mbox{with}\quad
   \sum_{k=1}^nL_k(t)<1}
for all $t\in H$ and $x_k,y_k\in B$ with continuous functions $L_k\colon X\to]0,+\infty[$. Then, there exists a unique
continuous solution $f\colon X\to B$ of the Presi\'c-type functional equation \eqref{Pres}.}

The Presi\'c equation was thoroughly treated and discussed in the monograph \cite{Cze80} by Czerwik. Our results are parallel to those
in \cite{Cze80}, however, our assumptions are less technical and our theorems are more general in many aspects, but these results are
not comparable in general.

\textbf{Acknowledgement.} We wish to express our gratitude to professor \textsc{Karol Baron}, who called our attention to Czerwik's monograph
and sent us its hard copy, and to professor \textsc{\'Arp\'ad Sz\'az}, who suggested to use the relation terminology for the formulation of
our main results.

%\bibliographystyle{amsplain}
%\bibliography{fixponthoz,additional}

\end{document}